\magnification=\magstep1
\input amstex
\documentstyle{amsppt}

\leftheadtext{Thomas Jech and Saharon Shelah}   

\hsize=5in
\vsize=7.5in

\def\pf{\hfill $\square$}
\def\cf{\text{cf}}
\def\c{\cite}
\def\ws{\widehat{S}}
\def\wt{\widehat{T}}
\def\doc{\dot{C}}
\def\dof{\dot{F}}
\def\dop{\dot{P}}
\def\dq{\dot{Q}}
\def\dor{\dot{R}}
\def\ds{\dot{S}}
\def\ova{\overline{a}}
\def\ovp{\overline{p}}
\def\cupn{\overset{\infty}\to{\underset{n=0}\to\bigcup}}
\def\limn{\underset{n}\to\lim}

\topmatter
\title
	On Reflection of Stationary Sets in $\Cal P_\kappa \lambda$
\endtitle
\author
	Thomas Jech and Saharon Shelah   \\
	\ \ \ \ \\
	{\rm Department of Mathematics  \\
	The Pennsylvania State University  \\
	University Park, PA 16802  \\
	\ \ \ \   \\
	Institute of Mathematics  \\
	The Hebrew University  \\
	Jerusalem, Israel}
	\ \ \ \ \ \\
	\ \ \ \ \ \\
\endauthor
\abstract
	Let $\kappa$ be an inaccessible cardinal, and let $E_0 =
	\{x \in \Cal P_\kappa \kappa^+ : \cf\; \lambda_x = \cf\; \kappa_x\}$ and
	$E_1 = \{x \in \Cal P_\kappa \kappa^+ :   \kappa_x$ is 
	regular and $\lambda_x = \kappa_x^+\}$.  It is consistent 
	that the set $E_1$ is stationary and that every 
	stationary subset of $E_0$ reflects at almost every 
	$a \in E_1$.
\endabstract
\endtopmatter

\document
\footnote"{}"{Supported by NSF grants DMS-9401275 and DMS 97-04477.}

\newpage

\baselineskip 20pt
\subhead
1. Introduction
\endsubhead

We study reflection properties of stationary sets in the space $\Cal P_\kappa
\lambda$ where $\kappa$ is an inaccessible cardinal.  Let $\kappa$ be a regular
uncountable cardinal, and let $A \supseteq \kappa$.  The set $\Cal P_\kappa A$
consists of all $x \subset A$ such that $|x| < \kappa$.  Following \c3, a
set $C\subseteq \Cal P_\kappa A$ is {\it closed unbounded} if it is
$\subseteq$-cofinal and closed under unions of chains of length $< \kappa$;
$S\subseteq \Cal P_\kappa A$ is {\it stationary} if it has nonempty intersection
with every closed unbounded set.  Closed unbounded sets generate a
normal $\kappa$-complete filter, and we use the phrase ``almost all $x$''
to mean all $x \in \Cal P_\kappa A$ except for a nonstationary set.

Almost all $x \in \Cal P_\kappa A$ have the property that $x \cap \kappa$ is an
ordinal.  Throughout this paper we consider only such $x$'s, and
denote $x \cap \kappa = \kappa_x$.  If $\kappa$ is inaccessible then for almost all
$x$, $\kappa_x$ is a limit cardinal (and we consider only such $x$'s.)  By
\c5, the closed unbounded filter on $\Cal P_\kappa A$ is generated by the sets
$C_F = \{x : x \cap \kappa \in \kappa$ and $F(x^{<\omega})\subseteq x\}$ where
$F$ ranges over functions $F : A^{<\omega}\to A$.  It follows that a
set $S \subseteq \Cal P_\kappa A$ is stationary if and only if every model $M$
with universe $\supseteq A$ has a submodel $N$ such that $|N| < \kappa$, $N
\cap \kappa \in \kappa$ and $N \cap A \in S$.  In most applications, $A$ is
identified with $|A|$, and so we consider $\Cal P_\kappa \lambda$ where
$\lambda$ is a cardinal, $\lambda > \kappa$.  For $x \in \Cal P_\kappa \lambda$ we
denote $\lambda_x$ the order type of $x$.

We are concerned with {\it reflection} of stationary sets.  Reflection
properties of stationary sets of ordinals have been extensively studied,
starting with \c7.  So have been reflection principles for stationary
sets in $\Cal P_{\omega_1} \lambda$, following \c2.  In this paper we
concentrate on $\Cal P_\kappa \lambda$ where $\kappa$ is inaccessible.

\definition
{Definition}  Let $\kappa$ be an inaccessible and let $a \in \Cal P_\kappa \lambda$
be such that $\kappa_a$ is a regular uncountable cardinal.  A stationary 
set $S \subseteq \Cal P_\kappa\lambda$ {\it reflects} at $a$ if the set $S
\cap \Cal P_{\kappa_a}a$ is a stationary set in $\Cal P_{\kappa_a} a$.
\enddefinition

The question underlying our investigation is to what extent can
stationary sets reflect.
There are some limitations associated with cofinalities.  For
instance, let $S$ and $T$ be stationary subsets of $\lambda$ such that
every $\alpha \in S$ has cofinality $\omega$, every $\gamma \in T$ has
cofinality $\omega_1$, and for each $\gamma \in T$, $S \cap \gamma$ is
a nonstationary subset of $\gamma$ (cf. \c4).  Let $\ws = \{ x \in \Cal P_\kappa
\lambda : \sup x \in S\}$ and $\wt = \{a \in \Cal P_\kappa \lambda : \sup a \in
T\}$. Then $\ws$ does not reflect at any $a\in \wt$.

Let us consider the case when $\lambda = \kappa^+$.  As the example
presented above indicates, reflection will generally fail when dealing
with the $x$'s for which $\cf\;\lambda_x < \kappa_x$, and so we restrict
ourselves to the (stationary) set
$$
	\{ x \in \Cal P_\kappa \lambda \;:\; \cf\; \kappa_x \leq \cf \;\lambda_x\}
$$
Since $\lambda = \kappa^+$, we have $\lambda_x \leq \kappa_x^+$ for almost
all $x$.

Let 
$$
\aligned
	E_0 & = \{ x \in \Cal P_\kappa \kappa^+ \;:\; 
\kappa_x \text{ is a limit cardinal and }
\cf\; \kappa_x = \cf\; \lambda_x\}\,,  \\
	E_1 & = \{ x \in \Cal P_\kappa \kappa^+ \;:\; 
\kappa_x \text{ is inaccessible and }\lambda_x=\kappa_x^+\}\,. 
\endaligned
$$
The set $E_0$ is stationary, and if $\kappa$ is a large cardinal
(e.g. $\kappa^+$-supercompact) then $E_1$ is stationary; the statement
``$E_1$ is stationary'' is itself a large cardinal property (cf. \c1).
Moreover, $E_0$ reflects at almost every $a \in E_1$ and
consequently, reflection of stationary subsets of $E_0$ at elements
of $E_1$ is a prototype of the phenomena we propose to investigate.

Below we prove the following theorem:

\proclaim
{1.2. \ Theorem}  Let $\kappa$ be a supercompact cardinal.  There is a generic
extension  in which
\roster
\item"{(a) }" the set $E_1 = \{x \in \Cal P_\kappa \kappa^+ : \kappa_x$ is 
inaccessible and $\lambda_x = \kappa_x^+\}$ is stationary, and
\item"{(b) }" for every stationary set $S \subseteq E_0$, the set
	$\{a \in E_1 : S \cap \Cal P_{\kappa_a} a$ is nonstationary in
	$\Cal P_{\kappa_a} a\}$ is nonstationary.
\endroster
\endproclaim

A large cardinal assumption in Theorem 1.2 is necessary.  As mentioned
above, (a) itself has large cardinal consequences.  Moreover, (b) 
implies reflection of stationary subsets of the set $\{ \alpha <
\kappa^+ : \cf\; \alpha < \kappa\}$, which is also known to be strong
(consistency-wise).

\subhead
2. \ Preliminaries
\endsubhead

We shall first state several results that we shall use in the proof of
Theorem 1.2.

We begin with a theorem of Laver that shows that supercompact cardinals
have a $\diamondsuit$-like property:

\proclaim
{2.1. \ Theorem}  \c6 If $\kappa$ is supercompact then there is a 
function $f: \kappa \to V_\kappa$ such that for every $x$ there exists an
elementary embedding $j : V \to M$ with critical point $\kappa$ such
that $j$ witnesses a prescribed degree of supercompactness and
$(j(f))(\kappa) = x$.
\endproclaim

We say that the function $f$ has {\it Laver's property}. 
\definition
{2.2. \ Definition}
A forcing notion is {\it $<\kappa$-strategically closed} if for every
condition $p$, player $I$ has a winning strategy in the following game
of length $\kappa$:  \ Players I and II take turns to play a descending
$\kappa$-sequence of conditions $p_0 > p_1 > \cdots > p_{\xi} > \cdots, 
\xi < \kappa$, with $p > p_0$, such that player I moves at limit stages. 
Player I wins if for each limit $\lambda < \kappa$, the sequence
$\{p_{\xi}\}_{\xi < \lambda}$ has a lower bound.
\enddefinition

It is well known that forcing with a $<\kappa$-strategically closed notion
of forcing does not add new sequences of length $<\kappa$, and that every
iteration, with $<\kappa$-support, of $<\kappa$-strategically closed forcing
notions is $<\kappa$-strategically closed.

\definition
{2.3. \ Definition} \c8 A forcing notion satisfies the 
$<\kappa${\it -strategic}-$\kappa^+${\it -chain condition} if for every limit 
ordinal $\lambda < \kappa$, player I has a winning strategy in
the following game of length $\lambda$:

Players I and II take turns to play, simultaneously for each $\alpha
< \kappa^+$ of cofinality $\kappa$, descending $\lambda$-sequences of conditions
$p_0^{\alpha} > p_1^{\alpha} > \cdots > p_{\xi}^{\alpha} > \cdots, \xi
< \lambda$, with player II moving first and player I moving at limit 
stages.  In addition, player I chooses, at stage $\xi$, a closed 
unbounded set $E_{\xi} \subset \kappa^+$ and a function $f_{\xi}$ such
that for each $\alpha < \kappa^+$ of cofinality $\kappa$, $f_{\xi}(\alpha) < 
\alpha$.

Player I wins if for each limit $\eta < \lambda$, each sequence
$\langle p_{\xi}^{\alpha} : \xi < \eta\rangle$ has a lower bound,
and if the following holds: \ for all $\alpha,\beta \in 
\bigcap_{\xi < \lambda} E_{\xi}$, if $f_{\xi}(\alpha)
= f_{\xi}(\beta)$ for all $\xi < \lambda$, then the sequences $\langle
p_{\xi}^{\alpha} : \xi < \lambda \rangle$ and $\langle p_{\xi}^{\beta}
: \xi < \lambda\rangle$ have a common lower bound.
\enddefinition

It is clear that property (2.3) implies the $\kappa^+$-chain condition.  
Every iteration with $<\kappa$-support, of
$<\kappa$-strategically $\kappa^+$-c.c. forcing notions satisfies the
$<\kappa$-strategic $\kappa^+$-chain condition. This is stated in \c8
and a detailed proof will appear in \c9.

In Lemmas 2.4 and 2.5 below, $H(\lambda)$ denotes the set of all sets
hereditarily of cardinality $< \lambda$.

\proclaim
{2.4. \ Lemma} Let $S$ be a stationary subset of $E_0$.  For every set
$u$ there exist a regular $\lambda > \kappa^+$, an elementary submodel $N$
of $\langle H(\lambda),\in, \Delta,u\rangle$ (where $\Delta$ is a well
ordering of $H(\lambda)$) such that $N \cap \kappa^+ \in S$, and a sequence
$\langle N_{\alpha} : \alpha < \delta\rangle$ of submodels of $N$ such
that $|N_{\alpha}|  < \kappa$ for every $\alpha$, $N \cap \kappa^+ =
\bigcup_{\alpha < \delta} (N_{\alpha} \cap \kappa^+)$ and for all
$\beta < \delta$, $\langle N_{\alpha} : \alpha < \beta\rangle \in N$.
\endproclaim 

\demo
{Proof} Let $\mu > \kappa^+$ be such that $u \in H(\mu)$, and let $\lambda
= (2^{\mu})^+$; let $\Delta$ be a well ordering of $H(\lambda)$.
There exists an elementary submodel $N$ of $\langle H(\lambda), \in ,
\Delta$) containing $u$, $S$ and $\langle H(\mu),\in , \Delta
\upharpoonright H(\mu)\rangle$ such that $N \cap \kappa^+ \in S$ and $N \cap \kappa$ 
is a strong limit cardinal; let $a = N \cap \kappa^+$.

Let $\delta = \cf\;\kappa_a$.  As $a\in S$, we have $\cf\;(\sup a) =
\delta$, and let $\gamma_{\alpha}$, $\alpha < \delta$, be an increasing
sequence of ordinals in $a - \kappa$, cofinal in $\sup a$.  Let $\langle
f_{\alpha} : \kappa \leq \alpha < \kappa^+\rangle \in N$ be such that each
$f_{\alpha}$ is a one-to-one function of $\alpha$ onto $\kappa$.  (Thus for
each $\alpha \in a$, $f_{\alpha}$ maps $a \cap \alpha$ onto $\kappa_a$.)
There exists an increasing sequence $\beta_{\alpha}$, $\alpha < \delta$,
of ordinals cofinal in $\kappa_a$, such that for each $\xi < \alpha$,
$f_{\gamma_{\alpha}} (\gamma_{\xi}) < \beta_{\alpha}$.

For each $\alpha < \delta$, let $N_{\alpha}$ be the Skolem hull
of $\beta_{\alpha} \cup \{\gamma_{\alpha}\}$ in $\langle H(\mu),
\in , \Delta \upharpoonright H(\mu), \langle f_{\alpha}\rangle\rangle$.  
$N_\alpha$
is an elementary submodel of $H(\mu)$ of cardinality $<\kappa_a$, and
$N_{\alpha} \in N$.  Also, if $\xi < \alpha$ then $\gamma_{\xi}
\in N_{\alpha}$ (because $f_{\gamma_{\alpha}} (\gamma_{\xi}) <
\beta_{\alpha}$) and so $N_{\xi} \subseteq N_{\alpha}$.
\enddemo

As $N \cap \kappa$ is a strong limit cardinal, it follows that for all
$\beta < \delta$, $\langle N_{\alpha} : \alpha < \beta\rangle \in N$.
Also, $N_{\alpha} \subseteq N$ for all $\alpha < \delta$, and it
remains to prove that $a \subseteq \bigcup_{\alpha < \delta} N_{\alpha}$.

As $\sup \{\beta_{\alpha} : \alpha < \delta\} = \kappa_a$, we have 
$\kappa_a \subseteq \bigcup_{\alpha < \delta} N_{\alpha}$.
If $\gamma \in a$, there exists a $\xi < \alpha < \delta$ such that
$\gamma < \gamma_{\xi}$ and $f_{\gamma_\xi}(\gamma) < \beta_{\alpha}$.
Then $\gamma_{\xi} \in N_{\alpha}$ and so $\gamma \in N_{\alpha}$.

\proclaim
{2.5. Lemma}  Let $S$ be a stationary subset of $E_0$ and let $P$
be a $<\kappa$-strategically closed notion of forcing.  Then $S$ remains
stationary in $V^P$. 
\endproclaim

\demo
{Proof} Let $\doc$ be a $P$-name for a club set in $\Cal P_\kappa \kappa^+$, and let
$p_0 \in P$.  We look for a $p \leq p_0$ that forces $S \cap \doc \neq
\emptyset$.

Let $\sigma$ be a winning strategy for I in the game (2.2).  By Lemma
2.4 there exist a regular $\lambda > \kappa^+$, an elementary submodel $N$
of $\langle H(\lambda), \epsilon,\Delta, P, p_0, \sigma, S,\doc\rangle$
(where $\Delta$ is a well-ordering) such that $|N| < \kappa$ and $N \cap \kappa^+
\in S$, and a sequence $\langle N_{\alpha} : \alpha < \delta\rangle$
of submodels of $N$ such that $|N_{\alpha}| < \kappa$ for every $\alpha$,
$N \cap \kappa^+ = \bigcup_{\alpha < \delta} (N_{\alpha} \cap
\kappa^+)$ and for all $\beta < \delta$, $\langle N_{\alpha} : \alpha <
\beta\rangle \in N$.

We construct a descending sequence of conditions $\langle p_{\alpha}
: \alpha < \delta\rangle$ below $p_0$ such that for all $\beta <
\delta$, $\langle p_{\alpha} : \alpha < \beta\rangle \in N$: at
each limit stage $\alpha$ we apply the strategy $\sigma$ to get
$p_{\alpha}$; at each $\alpha + 1$ let $q \leq p_{\alpha}$ be
the $\Delta$-least condition such that for some $M_{\alpha} \in
\Cal P_\kappa \kappa^+ \cap N$, $M_{\alpha} \supseteq N_{\alpha} \cap \kappa^+$, 
$M_{\alpha} \supseteq \bigcup_{\beta < \alpha} M_{\beta}$
and $q \Vdash M_{\alpha} \in \doc$ (and let $M_{\alpha}$ be the $\Delta$-least
such $M_\alpha$), and then apply $\sigma$ to get $p_{\alpha+1}$.
Since $M_\alpha \in N$, $N\models
|M_{\alpha}| < \kappa$ and so $M_{\alpha}\subseteq N$;
hence $M_{\alpha}\subseteq N \cap \kappa^+$.  Since for all $\beta <
\delta$, $\langle N_{\alpha} : \alpha < \beta \rangle \in N$,
the construction can be carried out inside $N$ so that for each
$\beta < \delta$, $\langle p_{\alpha} : \alpha < \beta\rangle \in N$.

As I wins the game, let $p$ be a lower bound for $\langle p_{\alpha} :
\alpha < \delta\rangle$; $p$ forces that $\doc \cap (N \cap \kappa^+)$ is
unbounded in $N \cap \kappa^+$ and hence $N \cap \kappa^+ \in \doc$.  Hence
$p \Vdash S \cap \doc \neq \emptyset$.  
\line{\hfill \pf}
\enddemo

\subhead
3. \ The forcing
\endsubhead

We shall now describe the forcing construction that yields Theorem 1.2.
Let $\kappa$ be a supercompact cardinal.

The forcing $P$ has two parts, $P = P_\kappa * \dop^\kappa$, where $P_\kappa$ is the
{\it preparation forcing} and $P^\kappa$ is the {\it main iteration}.  
The preparation forcing is an iteration of length $\kappa$, with Easton
support, defined as follows: \ Let $f : \kappa \to V_\kappa$ be a function
with Laver's property.  If $\gamma < \kappa$ and if $P_\kappa \upharpoonright
\gamma$ is the iteration up to $\gamma$, then the $\gamma^{\text{th}}$
iterand $\dq_{\gamma}$ is trivial unless $\gamma$ is inaccessible and
$f(\gamma)$ is a $P_\kappa \upharpoonright \gamma$-name for a
$<\gamma$-strategically closed forcing notion, in which case
$\dq_{\gamma} = f(\gamma)$ and $P_{\gamma + 1} = P_{\gamma} *
\dq_{\gamma}$.  Standard forcing arguments show that $\kappa$ remains
inaccessible in $V^{P_\kappa}$ and all cardinals and cofinalities above $\kappa$
are preserved.

The main iteration $\dop^\kappa$ is an iteration in $V^{P_\kappa}$, of length
$2^{(\kappa^+)}$, with $<\kappa$-support.  We will show that each iterand
$\dq_{\gamma}$ is $<\kappa$-strategically closed and satisfies the
$<\kappa$-strategic $\kappa^+$-chain condition.  This guarantees that $\dop^\kappa$ is
(in $V^{P_\kappa}$) $<\kappa$-strategically closed and satisfies the $\kappa^+$-chain
condition, therefore adds no bounded subsets of $\kappa$ and preserves all
cardinals and cofinalities.

Each iterand of $\dop^\kappa$ is a forcing notion $\dq_{\gamma} = Q(\ds)$
associated with a stationary set $\ds \subseteq \Cal P_\kappa \kappa^+$ in $V^{P_\kappa *
\dop_\kappa\upharpoonright\gamma}$, to be defined below.  By the usual
bookkeeping method we ensure that for every $P$-name $\ds$ for a
stationary set, some $\dq_{\gamma}$ is $Q(\ds)$.

Below we define the forcing notion $Q(S)$ for every stationary set
$S \subseteq E_0$; if $S$ is not a stationary subset of $E_0$ 
then $Q(S)$ is the trivial forcing. If $S$ is a stationary subset of $E_0$
then a generic for $Q(S)$ produces
a closed unbounded set $C \subseteq \Cal P_\kappa \kappa^+$ such that for every
$a \in E_1 \cap C$, $S \cap \Cal P_{\kappa_a} a$ is stationary in $\Cal P_{\kappa_a}a$.
Since $\dop^\kappa$ does not add bounded subsets of $\kappa$, the forcing
$Q(\ds)$ guarantees that in $V^P$, $\ds$ reflects at almost every
$a \in E_1$.  The crucial step in the proof will be to show that the
set $E_1$ remains stationary in $V^P$.

To define the forcing notion $Q(S)$ we use certain models with 
universe in $\Cal P_\kappa \kappa^+$.  We first specify what models we use:

\definition
{3.1. Definition}  A {\it model} is a structure $\langle M,\pi,\rho\rangle$
such that
\roster
\item"{(i) }"  $M \in \Cal P_\kappa \kappa^+$; $M \cap \kappa = \kappa_M$ 
is an ordinal and
	$\lambda_M =$ the order type of $M$ is at most $|\kappa_M|^+$
\item"{(ii) }" $\pi$ is a two-place function; $\pi(\alpha,\beta)$
	is defined for all $\alpha \in M-\kappa$ and $\beta \in M\cap \alpha$.
	For each $\alpha \in M - \kappa$, $\pi_{\alpha}$ is the function
	$\pi_{\alpha}(\beta) = \pi(\alpha,\beta)$ from $M\cap\alpha$ onto
$M\cap\alpha$, and moreover, $\pi_{\alpha}$ maps $\kappa_M$ onto $M \cap \alpha$.
\item"{(iii) }" $\rho$ is a two-place function; $\rho(\alpha,\beta)$
	is defined for all $\alpha \in M - \kappa$ and $\beta < \kappa_M$.
	For each $\alpha\in M-\kappa$, $\rho_{\alpha}$ is the function 
	$\rho_{\alpha}(\beta) = \rho(\alpha,\beta)$ from $\kappa_M$ into
$\kappa_M$, and $\beta 
	\leq \rho_{\alpha}(\beta) < \kappa_M$ for all $\beta < \kappa_M$.
\endroster 
Two models $\langle M, \pi^M,\rho^M\rangle$ and $\langle N,\pi^N,
\rho^N\rangle$ are {\it coherent} if $\pi^M(\alpha,\beta) = \pi^N
(\alpha,\beta)$ and $\rho^M(\alpha,\beta) = \rho^N(\alpha,\beta)$
for all $\alpha,\beta \in M \cap N$.  $M$ is a {\it submodel} of
$N$ if $M \subseteq N$, and $\pi^M \subseteq \pi^N$ and $\rho^M
\subseteq \rho^N$.
\enddefinition

\proclaim
{3.2.  Lemma}  Let $M$ and $N$ be coherent models with $\kappa_M \leq
\kappa_N$.  If $M \cap N$ is cofinal in $M$ (i.e. if for all $\alpha \in
M$ there is a $\gamma \in M \cap N$ such that $\alpha < \gamma$),
then $M \subseteq N.$ 
\endproclaim

\demo
{Proof} Let $\alpha \in M$; let $\gamma \in M \cap N$ be such that
$\alpha < \gamma$.  As $\pi_{\gamma}^M$ maps $\kappa_M$ onto $M \cap
\gamma$, there is a $\beta < \kappa_M$ such that $\pi_{\gamma}^M(\beta) =
\alpha$.  Since both $\beta$ and $\gamma$ are in $N$, we have $\alpha
= \pi^M (\gamma,\beta) = \pi^N(\gamma,\beta) \in N$.
\enddemo

We shall now define the forcing notion $Q(S)$:

\definition
{3.3  Definition}  Let $S$ be a stationary subset of the set
$E_0 = \{x \in \Cal P_\kappa \kappa^+ : \kappa_x$ is a limit cardinal and $\cf\;
\lambda_x = \cf\;\kappa_x\}$.  A {\it forcing condition} in $Q(S)$ 
is a model $M = \langle M,\pi^M, \rho^M\rangle$ such that
\roster
\item"{(i) }"  $M$ is $\omega$-closed, i.e. for every ordinal
	$\gamma$, if $\cf\;\gamma = \omega$ and $\sup(M\cap \gamma)
	=\gamma$, then $\gamma \in M$;
\item"{(ii) }" For every $\alpha \in M - \kappa$ and $\beta < \kappa_M $,
	if $\kappa_M \leq \gamma < \alpha$, and if $\{\beta_n : n <
	\omega\}$ is a countable subset of $\beta$ such that
	$\gamma = \sup \{\pi_{\alpha}^M (\beta_n):n < \omega\}$,
	then there is some $\zeta < \rho_{\alpha}^M(\beta)$ such that
	$\gamma = \pi_{\alpha}^M(\zeta)$.
\item"{(iii) }" For every submodel $a \subseteq M$, if
$$
	a \in E_1 = \{ x \in \Cal P_\kappa \kappa^+ : \kappa_x \text{ is inaccessible and }
	\lambda_x = \kappa_x^+\},
$$
	then $S \cap \Cal P_{\kappa_a}a$ is stationary in $\Cal P_{\kappa_a}a$.  
\endroster

A forcing condition $N$ is {\it stronger} than $M$ if
$M$ is a submodel of $N$ and $|M| < |\kappa_N|$.
\enddefinition

The following lemma guarantees that the generic for $Q_S$ is unbounded
in $\Cal P_\kappa \kappa^+$.

\proclaim
{3.4.  Lemma}  Let $M$ be a condition and let $\delta < \kappa$ and
$\kappa \leq \varepsilon < \kappa^+$.  Then there is a condition $N$ stronger
than $M$ such that $\delta \in N$ and $\varepsilon \in N$.
\endproclaim

\demo
{Proof}  Let $\lambda < \kappa$ be an inaccessible cardinal, such that
$\lambda \geq \delta$ and $\lambda > |M|$.  We let $N = M \cup
\lambda \cup \{\lambda\} \cup \{\varepsilon\}$; thus $\kappa_N = \lambda 
+ 1$, and $N$ is $\omega$-closed.  We extend $\pi^M$ and $\rho^M$
to $\pi^N$ and $\rho^N$ as follows:

If $\kappa \leq \alpha < \varepsilon$ and $\alpha \in M$, we let
$\pi_{\alpha}^N(\beta) = \beta$ for all $\beta \in N$ such that $\kappa_M
\leq \beta \leq \lambda$.  If $\alpha \in M$ and $\varepsilon < \alpha$, we
define $\pi_{\alpha}^N$ so that $\pi_{\alpha}^N$ maps $\kappa_N - \kappa_M$ onto
$(\kappa_N - \kappa_M) \cup \{\varepsilon\}$.  For $\alpha = \varepsilon$, we define
$\pi_{\varepsilon}^N$ in such a way that $\pi_{\varepsilon}^N$ maps
$\lambda$ onto $N \cap \varepsilon$.

Finally, if $\alpha,\beta \in N$, $\beta < \kappa \leq \alpha$, and if 
either $\alpha = \varepsilon$ or $\beta \geq \kappa_M$ we let $\rho_{\alpha}^N
(\beta) = \lambda$.

Clearly, $N$ is a model, $M$ is a submodel of $N$, and $|M| < |\kappa_N|$.
Let us verify (3.3.ii).  This holds if $\alpha \in M$, so let $\alpha
=\varepsilon$.  Let $\beta \leq \lambda$, let $\{\beta_n : n < \omega\}
\subseteq \beta$ and let $\gamma = \sup \{\pi_{\varepsilon}^N(\beta_n)
: n < \omega\}$ be such that $\kappa \leq \gamma < \varepsilon$.  There is
a $\zeta < \lambda = \rho_{\varepsilon}^N(\beta)$ such that 
$\pi_{\varepsilon}^N(\zeta) = \gamma$, and so (3.3.ii) holds.  

To complete the proof that $N$ is a forcing condition, we verify
(3.3. iii).  This we do by showing that if $a \in E_1$ is a submodel
of $N$ then $a \subseteq M$.

Assume that $a \in E_1$ is a submodel of $N$ but $a \nsubseteq M$.
Thus there are $\alpha,\beta \in a$, $\beta < \kappa \leq \alpha$ such that
either $\alpha = \varepsilon$ or $\beta \geq \kappa_M$.  Then $\rho_{\alpha}^a
(\beta ) = \rho_{\alpha}^N(\beta) = \lambda$ and so $\lambda \in a$,
and $\kappa_a = \lambda + 1$.  This contradicts the assumption that $\kappa_a$
is an inaccessible cardinal.  \pf 
\enddemo

Thus if $G$ is a generic for $Q_S$, let $\langle M_G,\pi_G,\rho_G
\rangle$ be the union of all conditions in $G$.  Then for every 
$a \in E_1$, that is a submodel of $M_G$, 
$S \cap \Cal P_{\kappa_a} a$ is stationary
in $\Cal P_{\kappa_a} a$.  Thus $Q_S$ forces that $S$ reflects at all but
nonstationary many $a \in E_1$.  

We will now prove that the forcing
$Q_S$ is $<\kappa$-strategically closed.  The key technical devices are the
two following lemmas.

\proclaim
{Lemma 3.5}  Let $M_0 > M_1 > \dots > M_n > \dots$ be an $\omega$-sequence
of conditions.  There exists a condition $M$ stronger than all the 
$M_n$, with the following property:
$$
\split
\text{If $N$ is any model coherent with $M$ such that there
exists some $\gamma \in N \cap M$} \\
\text{but $\gamma \notin\cupn M_n$, then $\kappa_N > \limn\; \kappa_{M_n}$.}
\endsplit\tag 3.6
$$
\endproclaim

\demo
{Proof} Let $A = \cupn M_n$ and $\delta = A \cap \kappa = \limn\; \kappa_{M_n}$,
and let $\pi^A = \cupn \pi^{M_n}$ and $\rho^A = \cupn \rho^{M_n}$.  We
let $M$ be the $\omega$-closure of $(\delta + \delta) \cup A$; hence
$\kappa_M = \delta + \delta + 1$.
To define $\pi^M$, we first define $\pi_{\alpha}^M \supset
\pi_{\alpha}^A$ for $\alpha\in A$ in such a way that $\pi_\alpha^M$ maps
$\delta+\delta$ onto $M\cap\alpha$. When $\alpha \in M - A$
and $\alpha \geq \kappa$, we have $|M \cap \alpha | = |\delta|$ and so
there exists a function $\pi_{\alpha}^M$ on $M\cap\alpha$ that
maps $\delta + \delta$ onto $M \cap \alpha$; we let $\pi_{\alpha}^M$
be such, with the additional requirement that $\pi_{\alpha}^M(0) =
\delta$.  To define $\rho^M$, we let $\rho^M \supset \rho^A$ be such that
$\rho^M(\alpha,\beta) = \delta + \delta$ whenever either $\alpha
\notin A$ or $\beta \notin A$.

We shall now verify that $M$ satisfies (3.3. ii).  Let $\alpha,\beta 
\in M$ be such that $\alpha \geq \kappa$ and $\beta < \kappa$ and let $\gamma
\in M$, $\kappa \leq \gamma < \alpha$, be an $\omega$-limit point of the
set $\{\pi_{\alpha}^M(\xi) : \xi < \beta\}$.  We want to show that
$\gamma = \pi_{\alpha}^M(\eta)$ for some $\eta <
\rho_{\alpha}^M(\beta)$.  If both $\alpha$ and $\beta$ are in $A$ then
this is true, because $\alpha,\beta \in M_n$ for some $n$, and $M_n$
satisfies (3.3 ii).  If either $\alpha \notin A$ or $\beta \notin A$
then $\rho_{\alpha}^M(\beta) = \delta + \delta$, and since $\pi_{\alpha}^M$
maps $\delta + \delta$ onto $M \cap \alpha$, we are done.

Next we verify that $M$ satisfies (3.6).  Let $N$ be any model
coherent with $M$, and let $\gamma \in M \cap N$ be such that $\gamma
\notin A$.  If $\gamma < \kappa$ then $\gamma \geq \delta$ and so $\kappa_N >
\delta$.  If $\gamma \geq \kappa$ then $\pi_{\gamma}^M(0) = \delta$,
and so $\delta = \pi_{\gamma}^N(0) \in N$, and again we have $\kappa_N
> \delta$.

Finally, we show that for every $a \in E_1$, if $a \subseteq M$ then
$S \cap \Cal P_{\kappa_a} a$ is stationary.  We do this by showing that for
every $a \in E_1$, if $a \subseteq M$ then $a \subseteq M_n$ for
some $M_n$.

Thus let $a \subseteq M$ be such that $\kappa_a$ is regular and
$\lambda_a = \kappa_a^+$.  As $\kappa_a \leq \kappa_M = \delta + \delta + 1$,
it follows that $\kappa_a < \delta$ and so $\kappa_a < \kappa_{M_{n_0}}$ for
some $n_0$.  Now by (3.6) we have $a \subseteq \cupn M_n$, and
since $\lambda_a$ is regular uncountable, there exists some 
$n\geq n_0$ such that $M_n \cap a$ is cofinal in $a$.
It follows from Lemma 3.2 that $a \subseteq M_n$.
\enddemo

\proclaim
{Lemma 3.7}
Let $\lambda < \kappa$ be a regular uncountable cardinal and let 
$M_0 > M_1 > \cdots > M_{\xi} > \cdots$, $\xi < \lambda$,  
be a $\lambda$-sequence of conditions with the property that
for every $\eta < \lambda$ of cofinality $\omega$,
$$
\split
\text{If $N$ is any model coherent with $M_{\eta}$ such that
there exists some $\gamma \in N \cap M_{\eta}$} \\
\text{but $\gamma \notin
\bigcup_{\xi < \eta} M_{\xi}$, then $\kappa_N > \underset{\xi \to
\eta}\to\lim\; \kappa_{M_{\xi}}$.}
\endsplit\tag 3.8
$$
Then $M = \bigcup_{\xi < \lambda} M_{\xi}$ is a condition.
\endproclaim 

\demo
{Proof} It is clear that $M$ satisfies all the requirements for a
condition, except perhaps (3.3 iii).  ($M$ is $\omega$-closed because
$\lambda$ is regular uncountable.)  Note that because $|M_{\xi}| <
\kappa_{M_{\xi + 1}}$ for all $\xi < \lambda$, we have $|M| = |\kappa_M|$.

We shall prove (3.3 iii) by showing that for every $a \in E_1$, if
$a \subseteq M$, then $a \subseteq M_{\xi}$ for some $\xi < \lambda$.
Thus let $a \subseteq M$ be such that $\kappa_a$ is regular and $\lambda_a
= \kappa_a^+$.

 As $\lambda_a = |a| \leq |M| = |\kappa_M|$, it follows that
$\kappa_a < \kappa_M$ and so $\kappa_a < \kappa_{M_{\xi_0}}$ for some $\xi_0 < \lambda$.
We shall prove that there exists some $\xi \geq \xi_0$ such that
$M_{\xi} \cap a$ is cofinal in $a$; then by Lemma 3.2, $a \subseteq
M_{\xi}$.

We prove this by contradiction.  Assume that no $M_{\xi} \cap a$ is
cofinal in $a$.  We construct sequences $\xi_0 < \xi_1 < \cdots <
\xi_n < \cdots$ and $\gamma_1 < \gamma_2 < \cdots < \gamma_n < \cdots$
such that for each $n$,
$$
	\gamma_n \in a\,,\qquad \gamma_n > \sup(M_{\xi_n} \cap a)\,,
	\qquad \text{ and } \quad\gamma_n \in M_{\xi_{n+1}}
$$
Let $\eta = \lim_n \xi_n$ and $\gamma = \lim_n \gamma_n$.  We claim
that $\gamma \in a$.  

As $\lambda_a$ is regular uncountable, there
exists an $\alpha \in a$ such that $\alpha > \gamma$.  Let $\beta_n$, $n
\in \omega$, be such that $\pi_{\alpha}^a (\beta_n) = \gamma_n$, and
let $\beta < \kappa_a$ be such that $\beta > \beta_n$ for all $n$.  As $M$
satisfies (3.3. ii), and $\gamma = \sup \{\pi_{\alpha}^M(\beta_n) : n
< \omega\}$, there is some $\zeta < \rho_{\alpha}^M(\beta)$ such that
$\gamma = \pi_{\alpha}^M(\zeta)$.  Since $\zeta < \rho_{\alpha}^M(\beta) =
\rho_{\alpha}^a(\beta) < \kappa_a$, we have $\zeta \in a$, and $\gamma =
\pi_{\alpha}^a(\zeta) \in a$.

Now since $\gamma \in a$ and $\gamma > \sup(M_{\xi_n}\cap a)$ 
we have $\gamma \notin M_{\xi_n}$, for all $n$.  As $M_{\eta}$ is
$\omega$-closed, and $\gamma_n \in M_{\eta}$ for each $n$, we have
$\gamma \in M_{\eta}$.  Thus by (3.8) it follows that $\kappa_a > \lim_n
\kappa_{M_{\xi_n}}$, a contradiction. \pf
\enddemo

\proclaim
{Lemma 3.9}  $Q_S$ is $<\kappa$-strategically closed.
\endproclaim

\demo
{Proof} In the game, player I moves at limit stages.  In order to win
the game, it suffices to choose at every limit ordinal $\eta$ of
cofinality $\omega$, a condition $M_{\eta}$ that satisfies (3.8).
This is possible by Lemma 3.5.  \pf
\enddemo

We shall now prove that $Q_S$ satisfies the $<\kappa$-strategic $\kappa^+$-chain
condition.  First a lemma:

\proclaim
{Lemma 3.10}  Let $\langle M_1, \pi_1, \rho_1\rangle$ and 
$\langle M_2, \pi_2, \rho_2\rangle$ be forcing conditions such 
that $\kappa_{M_1} = \kappa_{M_2}$  and that the models $M_1$ and $M_2$ 
are coherent.  Then the conditions are compatible.
\endproclaim

\demo
{Proof}  Let $\lambda < \kappa$ be an inaccessible cardinal such that
$\lambda > |M_1 \cup M_2|$ and let $M = M_1 \cup M_2 \cup \lambda
\cup \{\lambda\}$.  We shall extend $\pi_1 \cup \pi_2$ and $\rho_1
\cup \rho_2$ to $\pi^M$ and $\rho^M$ so that $\langle M,\pi^M,
\rho^M\rangle$ is a condition.

If $\alpha \in M_i - \kappa$, we define $\pi_{\alpha}^M \supset \pi_i$
so that $\pi_{\alpha}^M$ maps $\lambda - \kappa_{M_1}$ onto $M \cap 
\alpha$, and such that $\pi_{\alpha}^M(\beta) = \lambda$ whenever
$\kappa \leq \beta < \alpha$, $\alpha\in M_1-M_2$ and $\beta \in M_2 - M_1$
(or vice versa).  We define $\rho_{\alpha}^M
\supset \rho_i$ by $\rho_{\alpha}^M(\beta) = \lambda$ for $\kappa_{M_1}
\leq \beta \leq \lambda$.  It is easy to see that $M$ is an $\omega$-closed
model that satisfies (3.3 ii).

To verify (3.3 iii), we show that every $a \in E_1$ that is a
submodel of $M$ is either $a \subseteq M_1$ or $a \subseteq M_2$.
Thus let $a$ be a submodel of $M$, $a \in E_1$, such that 
neither $a \subseteq M_1$ nor $a \subseteq M_2$.  First assume
that $\kappa_a \leq \kappa_{M_1}$.  Then there are $\alpha,\beta \in a$ 
such that $\kappa \leq \beta < \alpha$ and $\alpha \in M_1-M_2$ while
$\beta \in M_2-M_1$ (or vice versa).  But then $\pi^a(\alpha,\beta) =
\pi^M(\alpha,\beta) = \lambda$ which implies $\lambda \in a$, or
$\kappa_a = \lambda + 1$, contradicting the inaccessibility of $\kappa_a$.

Thus assume that $\kappa_a > \kappa_{M_1}$.  Let $\alpha \in a$ be such
that $\alpha \geq \kappa$, and then we have $\rho^a(\alpha,\kappa_{M_1}) =
\rho^M(\alpha,\kappa_{M_1}) = \lambda$, giving again $\lambda \in a$,
a contradiction.  \pf
\enddemo

\proclaim
{Lemma 3.11}  $Q_S$ satisfies the $<\kappa$-strategic $\kappa^+$-chain
condition.
\endproclaim

\demo
{Proof} Let $\lambda$ be a limit ordinal $<\kappa$ and consider the game
(2.3) of length $\lambda$.  We may assume that $\cf\;\lambda > \omega$.
In the game, player I moves at limit stages, and the key to winning is
again to make right moves at limit stages of cofinality $\omega$.
Thus let $\eta$ be a limit ordinal $<\lambda$, and let
$\{M_{\xi}^{\alpha} : \alpha < \kappa^+,\; \cf\; \alpha = \kappa\}$ be the set of
conditions played at stage $\xi$.

By Lemma 3.5, player I can choose, for each $\alpha$, a condition
$M_{\eta}^{\alpha}$ stronger than each $M_{\xi}^{\alpha}$, $\xi <
\eta$, such that $M_{\eta}^{\alpha}$ satisfies (3.8).  Then let
$E_{\eta}$ be the closed unbounded subset of $\kappa^+$ 
$$
	E_{\eta} = \{\gamma < \kappa^+ : M_{\eta}^{\alpha} \subset
	\gamma \quad \text{ for all } \quad \alpha < \gamma\}\,,
$$
and let $f_{\eta}$ be the function $f_{\eta}(\alpha) = M_{\eta}^{\alpha}
\upharpoonright \alpha$, this being the restriction of the model
$M_{\eta}^{\alpha}$ to $\alpha$.

We claim that player I wins following this strategy:  By Lemma 3.7, player
I can make a legal move at every limit ordinal $\xi < \lambda$, and
for each $\alpha$ (of cofinality $\kappa$), $M^{\alpha} = 
\bigcup_{\xi < \lambda} M_{\xi}^{\alpha}$ is a condition.
Let $\alpha < \beta$ be ordinals of cofinality $\kappa$ in 
$\bigcap_{\xi < \lambda} E_{\xi}$ such that $f_{\xi}(\alpha)
= f_{\xi}(\beta)$ for all $\xi < \lambda$.  Then $M^{\alpha} \subset
\beta$ and $M^{\beta} \upharpoonright \beta = M^{\alpha} \upharpoonright
\alpha$, and because the functions $\pi$ and $\rho$ have the property
that $\pi(\gamma,\delta) < \gamma$ and $\rho(\gamma,\delta) < \gamma$
for every $\gamma$ and $\delta$, it follows that $M^{\alpha}$ and $M^{\beta}$
are coherent models with $\kappa_{M^{\alpha}} = \kappa_{M^{\beta}}$.  By 
Lemma 3.10, $M^{\alpha}$ and $M^{\beta}$ are compatible
conditions.  \pf
\enddemo

\subhead
4. \ Preservation of the set $E_1$
\endsubhead

We shall complete the proof by showing that the set
$$
	E_1 = \{ x \in \Cal P_\kappa \kappa^+ : \kappa_x \text{ is inaccessible and }
	\lambda_x = \kappa_x^+\}
$$
remains stationary after forcing with $P = P_\kappa * \dop^\kappa$.

Let us reformulate the problem as follows: Let us show, working in
$V^{P_\kappa}$, that for every condition $p \in \dop^\kappa$ and every
$\dop^\kappa$-name $\dof$ for an operation $\dof : (\kappa^+)^{<\omega} \to \kappa^+$
there exists a condition $\ovp \leq p$ and a set $x \in E_1$ such that
$\ovp$ forces that $x$ is closed under $\dof$.

As $\kappa$ is supercompact, there exists by the construction of $P_\kappa$ and
by Laver's Theorem 2.1, an elementary embedding $j:V \to M$ with
critical point $\kappa$ that witnesses that $\kappa$ is $\kappa^+$-supercompact
and such that the $\kappa^{\text{th}}$ iterand of the
iteration $j(P_\kappa)$ in $M$ is (the name for) the forcing $\dop^\kappa$.
The elementary embedding $j$ can be extended, by a standard argument,
to an elementary embedding $j : V^{P_\kappa} \to M^{j(P_\kappa)}$.  Since $j$
is elementary, we can achieve our stated goal by finding, in
$M^{j(P_\kappa)}$, a condition $\ovp \leq j(p)$ and a set $x \in j(E_1)$
such that $\ovp$ forces that $x$ is closed under $j(\dof)$.

The forcing $j(P_\kappa)$ decomposes into a three step iteration $j(P_\kappa) =
P_\kappa * \dop^\kappa * \dor$ where $\dor$ is, in $M^{P_\kappa * \dop^\kappa}$, a $<
j(\kappa)$-strategically closed forcing.

Let $G$ be an $M$-generic filter on $j(P_\kappa)$, such that $p \in G$.
The filter $G$ decomposes into $G = G_\kappa * H * K$ where $H$ and $K$
are generics on $\dop^\kappa$ and $\dor$ respectively, and $p \in H$.
We shall find $\ovp$ that extends not just $j(p)$ but each member
of $j'' H$ ($\ovp$ is a {\it master condition}).  That will
guarantee that when we let $x = j'' \Cal P_\kappa \kappa^+$ 
(which is in $j(E_1)$)  then $\ovp$ forces
that $x$ is closed under $j(\dof)$: this
is because $\ovp \Vdash j (\dof) \upharpoonright x = j'' F_H$, where
$F_H$ is the $H$-interpretation of $\dof$. 

We construct $\ovp$, a sequence $\langle p_{\xi} : \xi < j(2^{\kappa^+})
\rangle$, by induction.  When $\xi$ is not in the range of $j$, we
let $p_{\xi}$ be the trivial condition; that guarantees that the
support of $\ovp$ has size $< j(\kappa)$.  So let $\xi = j(\gamma)$ be
such that $\ovp \upharpoonright \xi$ has been constructed.

Let $M$ the model $\bigcup \{j(N) : N\in H_{\gamma}\}$ where $H_{\gamma}$
is the $\gamma^{\text{th}}$ coordinate of $H$.  The $\gamma^{\text{th}}$
iterand of $\dop^\kappa$ is the forcing $Q(S)$ where $S$ is a stationary
subset of $E_0$.  In order for $M$ to be a condition in $Q(j(S))$,
we have to verify that for every submodel $a \subseteq M$, if
$a \in j(E_1)$ then $j(S)$ reflects at $a$.

Let $a \in j(E_1)$ be a submodel of $M$.  If $\kappa_a < \kappa_M = \kappa$, then
$a = j'' \ova = j(\ova)$ for some $\ova \in E_1$, and $\ova$ is a
submodel of some $N \in H_{\gamma}$.  As $S$ reflects at $\ova$
it follows that $j(S)$ reflects at $a$.

If $\kappa_a = \kappa$, then $\lambda_a = \kappa^+$, and $a$ is necessarily cofinal
in the universe of $M$, which is $j'' \kappa^+$.  By Lemma 3.2, we have
$a = M$, and we have to show that $j(S)$ reflects at $j''\kappa^+$.  This
means that $j''S$ is stationary in $\Cal P_\kappa(j'' \kappa^+)$, or equivalently,
that $S$ is stationary in $\Cal P_\kappa \kappa^+$.

We need to verify that $S$ is a stationary set, in the model
$M^{j(P_\kappa) * j(\dop_\kappa)\upharpoonright j(\gamma)}$, while we know that
$S$ is stationary in the model $V^{P_\kappa * \dop^\kappa\upharpoonright
\gamma}$.  However, the former model is a forcing extension of 
the latter by a $<\kappa$-strategically closed forcing, and the result
follows by Lemma 2.5.

\enddocument